\documentclass[11pt]{amsart}
\usepackage{amssymb}
\usepackage{amsfonts}
\usepackage[T1]{fontenc}
\usepackage{pstricks}
\usepackage{pstricks-add}
\usepackage{pst-node}
\usepackage{pst-coil}
\usepackage{graphicx}

\theoremstyle{plain}
\newtheorem{theorem}{Theorem}

\newtheorem{corollary}[theorem]{Corollary}

\newtheorem{lemma}[theorem]{Lemma}
\newtheorem{proposition}[theorem]{Proposition}

\theoremstyle{remark}

\newtheorem*{acknowledgement}{Acknowledgements}

\theoremstyle{definition}
\newtheorem{remark}[theorem]{Remark}
\newtheorem{definition}[theorem]{Definition}

\numberwithin{equation}{section}
\numberwithin{theorem}{section}

\newcommand{\rintfrac}[2]{\genfrac{\lceil}{\rceil}{}{1}{#1}{#2}}

\newcommand{\rz}{\rho_{ab}}
\DeclareMathOperator{\imz}{Im}

\DeclareMathOperator{\supp}{supp}
\newcommand{\ctg}{\cot}
\newcommand{\TLS}{Tristram--Levine signature}
\newcommand{\dl}{\Delta_{IJ}(x)}
\newcommand{\dlh}{\widehat{\Delta}_{IJ}(t)}
\newcommand{\suone}{\sum\limits_{\substack{\alpha<1\\\alpha\in\Sigma_{p,q}}}}
\newcommand{\suhalf}{\sum\limits_{\substack{\alpha<1/2\\\alpha\in\Sigma_{p,q}}}}
\newcommand{\suhalfone}{\sum\limits_{\substack{\alpha\in(1/2,1)\\\alpha\in\Sigma_{p,q}}}}
\newcommand{\rh}[1]{$\rho$\nobreakdash--\hspace{0pt}#1}
\newcommand{\chf}[2]{\chi_{\left(#1,#2\right)}}

\begin{document}
\title{A \rh invariant of iterated torus knots}
\author{Maciej Borodzik}
\address{Institute of Mathematics, University of Warsaw, ul. Banacha 2,
02-097 Warsaw, Poland}
\email{mcboro@mimuw.edu.pl}
\date{29 June 2009}
\thanks{The author is partially supported by the Foundation for Polish Science}
\begin{abstract}
We compute \rh invariant for iterated torus knots $K$ for the standard representation
$\pi_1(S^3\setminus K)\to\mathbb{Z}$ given by abelianisation. For algebraic knots, this invariant
 turns out to be very closely related to an invariant of a plane curve singularity, coming from algebraic geometry.
\end{abstract}
\maketitle

\section{Introduction}
A von Neumann \rh invariant (also called $L^2$--signature, or $L^2$--eta invariant) of a real closed 3--manifold $M$ is a real number $\rho_\phi(M)$ 
associated to every representation $\phi:\pi_1(M)\to\Gamma$,  where $\Gamma$ is any group satisfying
PTFA condition (see \cite[Definition~2.1]{COT1}). As a special case, if $K$ is a knot in a 3--sphere, and 
we consider representations of the fundamental group of the manifold $S^3_0(K)$ (i.e. a zero surgery on $K$),
then we can talk about the \rh{invariants} of knots. In particular, the representation 
$ab:\nobreak\pi_1(S^3\setminus K)\to\mathbb{Z}$, given by abelianization, 
gives rise to the representation $\tilde{ab}:\pi_1(S^3_0(K))\to\mathbb{Z}$
and the corresponding invariant, $\rz(K)$, turns out to be the integral over normalised unit circle 
of the \TLS{} of a knot.

The \rh invariants for knots have been introduced first in \cite{ChG}. They were then deeply studied in \cite{COT1}. 
In their seminal paper, the authors 
observed that they are a very subtle obstruction for some knots to be slice. Namely, let us be given a knot $K$ bounding a disk $D$ in the ball $B^4$.
Let $Y=\partial(B^4\setminus \nu(D))$, where $\nu$ denotes the tubular neighbourhood. Then $Y$ is canonically isomorphic to $S^3_0(K)$, and, for any 
representation $\phi:\pi_1(Y)\to\Gamma$ that can be extended to $\tilde{\phi}:\pi_1(B^4\setminus\nu(D))\to\Gamma$, 
the corresponding
\rh invariant must vanish. This allows to construct examples of non-slice knots, undistinguishable from slice knots
by previously known methods as the \TLS{} or the Casson--Gordon invariants.

The difficulty of computability of \rh invariants is the cost of their subtlety. Only in the first nontrivial
case of the representation given by $ab$,  there is a general
method of computing this invariant, namely integrating the \TLS. In papers \cite{COT2}, \cite{Ha}, and others,
these invariants were computed also for some other representations of the knot group. But there, the choice of knots
is very specific.

In this paper we focus on $\rz$--invariant and compute it for all iterated torus knots. The computation consists of integrating
the \TLS, which is not a completely trivial task. In fact, we do even more: we compute the Fourier transform
of the \TLS{} function of iterated torus knot. This transform can be expressed by a surprisingly simple formula. 
In particular, this method can be used to detect
knots, which are connected sums of iterated torus knots and which have identical \TLS.

What we find most interesting and striking about $\rz$ of algebraic knots,
is its relation with deep algebro-geometrical invariants of 
the plane curve singularity. We state this relation, in terms of a uniform bound (see Proposition~\ref{bound})
but, honestly speaking, we are far from understanding it. Moreover, this
relation is not that clear for algebraic links, as we show on an example.

The structure of the paper is the following. In Section~2 we recall, how to compute the \TLS{} for iterated torus knots and formulate
Theorem~\ref{core}. Then we deduce some of its corollaries. In Section~3 we prove Theorem~\ref{core}. In Section~4 we recall definitions
of some invariants of plane curve singularities and compare them to $\rz$ for algebraic knots. We end this section by computing the $\rz$
for a $(d,d)$ torus link, i.e., the link of singularity $x^d-y^d=0$.

We apologise the reader for not giving a definition of the \rh invariant. A precise definition from scratch, including 
necessary definitions
of twisted signature of a $4-$manifold, would make this paper twice as long. Instead we refer to \cite[Section~5]{COT1}, or, for
more detailed treatment, to a book by L\"uck \cite{Lu}.

We end this introduction by remarking that the $\rho$ invariants were also studied in the context of mixed Hodge structures
of hypersurface singularities. The $\eta$ invariant, defined, for instance, in \cite[Section I]{Ne1}, 
is closely related to  the $\rz$ invariant in the case of plane curve singularities. 
We refer to \cite{Ne2,Ne3} for the detailed study of this invariant.
\section{Tristram--Levine signature of torus knots}
We begin this section with some definitions, which we give also to fix the notation used in the article.
\begin{definition}\label{iterdef}
A knot is called an \emph{iterated torus knot} if it arises from an unknot by finitely many cabling operations. An
iterated torus knot is of type $(p_1,q_1,\dots,p_n,q_n)$ if it is a $(p_1,q_1)$ cable of $(p_2,q_2)$ cable
of \dots of $(p_n,q_n)$ cable of an unknot. Fore example, a torus knot $T_{p,q}$ is an iterated torus knot of type $(p,q)$.
\end{definition}
\begin{definition}
Let $K$ be a knot, $S$ its Seifert matrix. Let $\zeta\in\mathbb{C}$, $|\zeta|=1$. The \emph{\TLS}, $\sigma_K(\zeta)$ is the signature
of the hermitian form given by
\begin{equation}\label{eq:TLSig}
(1-\zeta)S+(1-\bar\zeta)S^T.
\end{equation}
\end{definition}
It is well-known that the form \eqref{eq:TLSig} is degenerate (i.e. has non-trivial kernel) if and only if
$\zeta$ is a root of the Alexander polynomial $\Delta_K$ of $K$. 
The function $\zeta\to \sigma_K(\zeta)$ is piecewise constant with possible jumps only at the roots of the Alexander polynomial
$\Delta_K(\zeta)$. The value of $\sigma_K$ at such root can \emph{a priori} be different then left or right limit of $\sigma_K$ at that
point. However, there are only finitely many such values and they do not influence the integral. As we do not want to take care
of this values, we introduce a very handy notion.
\begin{definition}
We shall say that two piecewise-constant functions from a unit circle (or a unit interval) to real numbers are \emph{almost equal} if they are equal at 
all but finitely
many points.
\end{definition}
We would like to compute $\rz$ for an iterated torus knot. We will use Proposition~5.1 from \cite{COT2},
which we can formulate as follows.
\begin{proposition}
For any knot $K\subset S^3$ we have
\[\rz(K)=\int_0^1\sigma_K(e^{2\pi ix})dx.\]
\end{proposition}
Therefore, what we have to do, is to compute the integral of the \TLS{} for an iterated torus knot. We begin
with recalling results from \cite{Li}, where the function $\sigma_K$ is computed for iterated torus knots.

Let $p,q$ be coprime positive integers. Let $x$ be in the interval $[0,1]$. Consider the set
\[\Sigma=\Sigma_{p,q}=\left\{\frac{k}{p}+\frac{l}{q}\colon 1\le k<p,\,1\le l<q\right\}\subset [0,2]\cap\mathbb{Q}.\]
The function $s_{p,q}(x)$ is defined as
\[s_{p,q}(x)=-2\#\Sigma\cap (x,x+1)+\#\Sigma.\]
\begin{lemma}[\cite{Li}]
If $\zeta=e^{2\pi ix}$ is not a root of the polynomial $(t^{pq}-1)(t-1)/(t^p-1)(t^q-1)$, then the \TLS{} of the torus knot
$T_{p,q}$ at $\zeta$ is equal to $s_{p,q}(x)$.
\end{lemma}
Therefore, computing the \rh invariant of a torus knot boils down to computing the integral of the function $s_{p,q}(x)$.
Before we do this,
let us show, how one can compute the Tristram--Levine signatures of an iterated torus knot. We shall need another lemma from \cite{Li}.
\begin{lemma}\label{cable}
Let $K$ be a knot and $K_{p,q}$ be the $(p,q)-$cable on $K$. Then for any $\zeta\in\mathbb{C}$, $|\zeta|=1$, we have
\[\sigma_{K_{p,q}}(\zeta)=\sigma_{K}(\zeta^q)+\sigma_{T_{p,q}}(\zeta).\]
\end{lemma}
This allows a recursive computation for an iterated torus knot. Namely, let for $r>1$
\[s_{p,q;r}(x)=s_{p,q}(\lfloor rx\rfloor).\]
\begin{corollary}\label{iterated}
Let $K$ be an iterated torus knot of type $(p_1,q_1,\dots,p_n,q_n)$. Let $x\in[0,1]$ be such that $e^{2\pi ix}$ is not a root
of the Alexander polynomial of $K$. Denote by $r_k=q_1\dots q_{k-1}$. Then
\[
\sigma_{K}(e^{2\pi i x})=\sum_{k=1}^ns_{p_k,q_k,r_k}(x).
\]
\end{corollary}
The core of this section is
\begin{theorem}\label{core}
For any $\beta\in\mathbb{C}$ which is not an integer divisible by $r$ we have
\begin{equation}\label{eq:formula1}
\int_0^1e^{\pi i\beta x}s_{p,q,r}(x)\,dx=\frac{2e^{\pi i\beta/2}\sin\frac{\pi\beta}{2}}{\pi\beta}n_{p,q;r}(\frac{\pi\beta}{2}),
\end{equation}
where
\[n_{p,q;r}(t)=\ctg\frac{t}{pqr}\ctg\frac{t}{r}-\ctg\frac{t}{pr}\ctg\frac{t}{qr}.\]
\end{theorem}
In particular, by taking a limit $\beta\to 0$ we get
\[
\int_0^1s_{p,q,r}=-\frac13(p-\frac1p)(q-\frac1q).
\]
\begin{remark}
The function $n_{p,q;r}(t)$ will be called \emph{normalised Fourier transform}.
\end{remark}
We prove Theorem~\ref{core} in Section~\ref{coreproof}. Now we pass to corollaries.
\begin{corollary}\label{corcore}
The $\rz$ invariant of an iterated torus knot is equal to
\[
-\frac13\sum_{k=1}^n(p_k-\frac1{p_k})(q_k-\frac{1}{q_k}).
\]
\end{corollary}
Apart of this corollary, 
Theorem~\ref{core} has its interest of its own. In fact, it might help to study possible cobordism relations between iterated torus knot. For
example, Litherland showed in \cite{Li}, that the connected sum of knots $T_{2,3}$, $T_{3,5}$ and a $(2,5)$-cable on $T_{2,3}$ 
has the same \TLS{} as
a $T_{6,5}$. It might be possible that  normalised Fourier transforms of torus knots can help studying similar phenomena.
This could be done as follows.
\begin{lemma}
Let us be given two finite sets $I$ and $J$ of triples of integers $\{p,q,r\}$. Then the difference
\begin{equation}\label{eq:int}
\dl:=\sum_{i\in I}s_{p_i,q_i;r_i}(x)-\sum_{j\in J}s_{p_j,q_j;r_j}(x)
\end{equation}
is almost equal to zero for $x\in[0,1]$, if and only if the difference
\begin{equation}\label{eq:trans}
\dlh:=\sum_{i\in I}n_{p_i,q_i;r_i}(t)-\sum_{j\in J}n_{p_j,q_j;r_j}(t)
\end{equation}
is equal to zero on some open subset in $\mathbb{C}$.
\end{lemma}
\begin{proof}[Sketch of proof]
The 'only if' part is trivial. To prove the 'if' part we observe that $\dlh\cdot\frac{t}{e^t\sin t}$ is, up to a multiplicative constant,
and up to rescaling of the parameter $t$, the Fourier transform of $\dl$, when we extend $\dl$ by $0$ to the whole real line. On the other
hand, vanishing of $\dlh$ on some open subset of $\mathbb{C}$ implies that it is everywhere $0$.
\end{proof}
\begin{proposition}\label{cornew}
The condition that $\dl$ is almost equal to zero is equivalent to the fact, that two following conditions are satisfied at once
\begin{itemize}
\item[(a)] $\sum_{i\in I}(p_i-\frac{1}{p_i})(q_i-\frac{1}{q_i})=\sum_{j\in J}(p_j-\frac{1}{p_j})(q_j-\frac{1}{q_j})$.
\item[(b)] For any $t_0$ such that $\pi r_kt_0\in\mathbb{Z}$ for some $k\in I\cup J$ the residuum at $t_0$ of $\dlh$ is zero.
\end{itemize}
\end{proposition}
\begin{remark}
If $T$ is the least common multiplier of $p_kq_kr_k$ for $k\in I\cup J$, then $T\pi$ is the period of $\dlh$. It follows that
the condition (b) involves only finitely many equations.
\end{remark}
\begin{proof}[Proof of Proposition~\ref{cornew}]
Vanishing of $\dlh$ clearly implies (b). The equality in (a) is equivalent to $3\widehat{\Delta}_{IJ}(0)=0$. We shall prove
that (b) implies that $\dlh$ is bounded on $\mathbb{C}$. This is done as follows.

Observe that, in general, $\dlh$ can have poles only at such $t_0$'s, that $\pi r_kt_0\in\mathbb{Z}$, for some $k\in I\cup J$.
Moreover, these poles are at most of order $1$: in fact, it is a matter of simple computation, that $n_{p,q;r}$ does not have
a pole of order $2$. Therefore, condition (b) implies that the $\dlh$ extends holomorphically across points $\frac{n}{\pi r_k}$,
where $k\in I\cup J$ and $n\in\mathbb{Z}$. As this function is periodic with real period, for any $\delta>0$ it is bounded
on the strip $|\imz t|\le\delta$ by some constants, depending of course of $\delta$.

A uniform bound on $\dlh$ for $|\imz t|\ge\delta$ results from the standard estimate $|\ctg t|^2\le 1+\frac{1}{(\imz t)^2}$. Hence, 
if (b) holds, then the function $\dlh$ 
is a bounded holomorphic function, by Liouville's theorem it is then constant. The condition~(a) implies then that it vanishes
at $0$, so it is zero everywhere.
\end{proof} 
\section{Proof of Theorem~\ref{core}}\label{coreproof}
To make computations at least a bit more transparent, let us first assume that $r=1$. 
The function $s_{p,q}$ can be expressed as the sum
\[s_{p,q}(x)=2\suhalf\chf{\alpha}{1-\alpha}(x)-2\suhalfone\chf{1-\alpha}{\alpha}(x),\]
where $\chf{a}{b}$ is the characteristic function of the interval $(a,b)$.
Therefore
\begin{equation}\label{eq:whatwewant}
\int_0^1s_{p,q}(x)e^{\pi i\beta x}dx=-\frac{2}{\pi i \beta}\suone e^{\pi i\alpha\beta}-e^{\pi i\beta(1-\alpha)}.
\end{equation}
We have
\[\suone e^{\pi i\alpha\beta}=\sum_{k=1}^{p-1}\sum_{\substack{l=1\\ l<q(1-k/p)}}^{q-1}e^{\pi i\beta(k/p+l/q)}.\]
The internal sum on the right hand side 
is the sum of geometric series (here we use the assumption that $\beta$ is not an integer) and can be expressed as
\[
\frac{1}{1-e^{\pi i\beta/q}}(e^{\pi i\beta k/p}-e^{\pi i\beta(k/p+l_k/q)}),
\]
where $l_k$ satisfies 
\[k/p+l_k/q>1>k/p+(l_k-1)/q.\]
So we have
\begin{equation}\label{eq:twosums}
\suone e^{\pi i\alpha\beta}=\frac{\sum\limits_{k=1}^{p-1}e^{\pi i\beta k/p}-\sum\limits_{k=1}^{p-1}e^{\pi i\beta(k/p+l_k/q)}}{1-e^{\pi i\beta/q}}.
\end{equation}
The first sum in the denominator is again a geometric series. As to the second one, let us denote 
\[\gamma_k=k/p+l_k/q.\]
Then $\gamma_k$'s have the following obvious properties
\begin{itemize}
\item[(a)] $\gamma_k$'s are all different;
\item[(b)] $1+\frac{1}{pq}\le\gamma_k\le 1+\frac{p-1}{pq}$;
\item[(c)] each $\gamma_k$ is of the form $1+a_k/pq$ with $a_k$ an integer.
\end{itemize}
By the Dirichlet principle the set $\{\gamma_1,\dots,\gamma_{p-1}\}$ is the same as the set $\{1+1/pq,\dots,1+(p-1)/pq\}$. 
Therefore, the second
sum in the denominator \eqref{eq:twosums}, upon reordering, can be expressed as
\[\sum_{m=1}^{p-1}e^{\pi i\beta (1+m/pq)},\]
which again is a geometric series. Putting things all together we get
\begin{equation}\label{eq:universalsum}
\suone e^{\pi i\alpha\beta}=\frac{1}{1-e^{\pi i\beta/q}}\left(\frac{e^{\pi i\beta/p}-e^{\pi i\beta}}{1-e^{\pi i\beta/p}}-
\frac{e^{\pi i\beta(1+1/pq)}-e^{\pi i\beta(1+1/q)}}{1-e^{\pi i\beta/pq}}\right).
\end{equation}
On the other hand, we have
\[
\suone e^{\pi i(1-\alpha)\beta}=e^{\pi i \beta}\suone e^{\pi i\alpha(-\beta)},
\]
and the sum on the right hand side is just \eqref{eq:universalsum} with $-\beta$ substituted in place of $\beta$.
Substituting this into \eqref{eq:whatwewant}, and applying the formula 
$e^{\pi i a}-e^{\pi ib}=2ie^{\pi i(a+b)/2}\sin\frac{\pi(a-b)}{2}$ several times, we arrive finally at
\[
\int_0^1 s_{p,q}(x)e^{\pi i\beta}dx=\frac{2e^{\pi i\beta/2}\sin{\frac{\pi\beta}{2}}}{\pi\beta}%
(\ctg\frac{\pi\beta}{2pq}\ctg\frac{\pi\beta}{2}-\ctg\frac{\pi\beta}{2p}\ctg\frac{\pi\beta}{2q}).
\]
To conclude the proof in the case $r>1$ we observe that
\begin{align*} 
s_{p,q;r}(x)=&2\suhalf \sum_{k=0}^{r-1}\chf{\frac{\alpha+k}{r}}{\frac{1-\alpha+k}{r}}(x)+\\
-&2\suhalfone \sum_{k=0}^{r-1}\chf{\frac{1-\alpha+k}{r}}{\frac{\alpha+k}{r}}(x)
\end{align*}
Thus
\begin{equation}\label{eq:involvedsum}
\int_0^1s_{p,q;r}e^{\pi i\beta x}=\frac{-2}{\pi i\beta}\suone \sum_{k=0}^{r-1}e^{\pi i\beta(\alpha/r+k/r)}-e^{\pi i\beta(1-\alpha/r-k/r)}.
\end{equation}
Now, for fixed $\alpha$ we have
\[\sum_{k=0}^{r-1}e^{\pi i\beta(\alpha/r+k/r)}=e^{\pi i\alpha(\beta/r)}\sum_{k=0}^{r-1}e^{\pi i\beta k/r}=e^{\pi i\alpha(\beta/r)}\frac{1-e^{\pi i\beta}}{1-e^{\pi i\beta/r}}.\]
Therefore, returning to \eqref{eq:involvedsum} we get
\[
\suone \sum_{k=0}^{r-1}e^{\pi i\beta(\alpha/r+k/r)}=\frac{1-e^{\pi i\beta}}{1-e^{\pi i\beta/r}}\suone e^{\pi i\alpha(\beta/r)}.
\]
We can use \eqref{eq:universalsum} again, substituting $\beta/r$ in place of $\beta$. Similarly we can deal with a sum of terms $e^{\pi i\beta(1-\alpha/r-k/r)}$.
Now straightforward but long computations yield the formula \eqref{eq:formula1}.
\section{Relation with algebraic invariants}
The setup in this section is the following. Let $(C,0)\subset\mathbb{C}^2$ be germ of a plane curve singularity with one branch.
This means that there exists a local parametrisation $C=(x(t),y(t))$, with $x$ and $y$ analytic functions in one variable with 
$x(0)=y(0)=0$. 
Let us assume that the Puiseux expansion of $y$ in fractional powers
of $x$ written is the multiplicative form (see \cite[page 49]{EN}) is
\[
y=x^{q_1/p_1}(a_1+x^{q_2/p_1p_2}(a_2+\ldots+x^{q_s/p_1p_2p_3\dots p_s}(a_s+\ldots))),
\] 
with $q_1>p_1$ (otherwise we switch $x$ with $y$), $\gcd(q_i,p_i)=1$ and $p_i,q_i>0$. 
The pairs $(p_1,q_1),\dots,(p_n,q_n)$ 
are called characteristic pairs (or Newton pairs) of the singularity. They completely determine the topological type of the singular point.
\begin{lemma}[see e.g. \cite{EN}]
Put $a_1=q_1$ and $a_{k+1}=p_{k+1}p_ka_k+q_{k+1}$. Then the link of the singularity $(C,0)$ is an iterated torus knot. More precisely, it is
a $(p_n,a_n)$ cable on $(p_{n-1},a_{n-1})$ cable on $\dots$ on $(p_1,a_1)$ torus knot 
\end{lemma}
\begin{remark}
The ordering of cables in \cite{EN} is different than in \cite{Li}. According to Definition~\ref{iterdef}, the link of the singuarity $(C,0)$ above
would be an iterated torus knot of type $(p_n,q_n,p_{n-1},q_{n-1},\dots,p_1,q_1)$.
\end{remark}
\begin{corollary}
The $\rz$ invariant of an algebraic knot is equal to
\begin{equation}\label{eq:tlint}
\rz=-\frac13\sum_{k=1}^n\left(a_kp_k-\frac{a_k}{p_k}-\frac{p_k}{a_k}+\frac{1}{p_ka_k}\right).
\end{equation}
\end{corollary}
It is on purpose that we wrote formula~\eqref{eq:tlint} in a different shape that in Corollary~\ref{corcore}.

Let us now resolve the above singularity. This means that we have a map $\pi:(X,E)\to (U,0)$, where $U$ is a neighbourhood
of $0$ in $\mathbb{C}^2$, $E$ is the exceptional divisor 
and $X$ is a complex surface. We require the strict transform $C'$ to
be smooth, $C'\cup E$ to have only normal crossings as singularities and the resolution 
to be minimal, so that we cannot blow-down any exceptional curve without
violating one of the two above assumptions.

Put $K=K_X$ the canonical divisor on $X$ and let $D=C'+E_{red}$. Here, the 
subscript 'red' means that we take a reduced divisor, i.e. coefficients with all components are equal
to $1$.
\begin{lemma}[\cite{OZ}]
Using the notation from this section, we have
\begin{equation}\label{eq:Mbar}
2\mu+(K+D)^2=a_1p_1-\rintfrac{a_1}{p_1}-\rintfrac{p_1}{a_1}+\sum_{k=2}^n\left(a_kp_k-\rintfrac{a_k}{p_k}\right),
\end{equation}
where $(K+D)^2$ denotes the self-intersection of the divisor $K+D$, and $\lceil x\rceil=\min(n\in\mathbb{Z},\,n\ge x)$.
\end{lemma}
On the one hand $(K+D)^2$ has a very natural meaning. Namely, at least for unibranched singularities,
this is the difference between the so called $\bar{M}$ number of singularity and the Milnor number $\mu$.
The $\bar{M}$ number, 
introduced in \cite{Or} and studied in \cite{BZ}, can be interpreted as a parametric codimension of a singular point, i.e. the number of locally
independent conditions, which are imposed on a curve given in parametric form, by the appearance of the singularity of given topological type.

On the other hand there is an apparent similarity of left hand sides of formulae \eqref{eq:tlint} and \eqref{eq:Mbar}. To make it even more similar, let us take
a Zariski--Fujita \cite{Fuj} decomposition of the divisor $K+D$. We have then
\[K+D=H+N\]
with $H$ nef (its intersection with any algebraic curve in $X$ is non-negative), $N$ effective and $N^2<0$, $H\cdot N'=0$ for any $N'$ supported on $\supp N$.
\begin{lemma}[\cite{OZ}]
\begin{equation}\label{eq:H}
2\mu+H^2=a_1p_1-\frac{a_1}{p_1}-\frac{p_1}{a_1}+\sum_{k=2}^n\left(a_kp_k-\frac{a_k}{p_k}\right).
\end{equation}
\end{lemma}
In the case of unibranched singularity, the quantity $H^2$ is the sum of Milnor number and so called $M$-number (without a bar) of singular point. Its importance
lies in the fact that the sum of $M$-numbers of all singular points of an algebraic curve in $\mathbb{C}P^2$ can be bounded from above by global topological data
of the curve, as genus and first Betti number (see \cite{BZ}). These bounds involve very deep Bogomolov--Miyaoka--Yau inequality from algebraic geometry.

Thus the following result seem to be a very mysterious and shows a deep link between knot theory and algebraic geometry.
\begin{proposition}\label{bound}
Let $\rz$ be the integral of the \TLS{} of an algebraic knot (see \eqref{eq:tlint}) and $H^2$ be like in \eqref{eq:H}. Then
\[0<-3\rz-(2\mu+H^2)<\frac29.\]
\end{proposition}
\begin{proof}
It easy to observe that
\[
\Delta:=-3\rz-(2\mu+H^2)=\frac{1}{a_1p_1}+\sum_{k=2}^n\left(\frac{1}{a_kp_k}-\frac{p_k}{a_k}\right).\]
On the one hand
\[\Delta\le\sum_{k=1}^n\frac{1}{a_kp_k}.\]
Recall that $a_{k+1}=a_kp_{k+1}p_k+q_{k+1}$, so $a_{k+1}p_{k+1}>a_kp_kp_{k+1}^2\ge 4a_kp_k$. Hence
\[\Delta\le\frac{1}{a_1p_1}\sum_{k=0}^{n-1}\frac{1}{4^k}<\frac{4}{3a_1p_1}.\]
But $a_1p_1\ge 6$, so one inequality is proved.

To prove in the second one, let us reorganise terms of $\Delta$ as follows
\[
\Delta=\sum_{k=1}^{n-1}\left(\frac{1}{a_kp_k}-\frac{p_{k+1}}{a_{k+1}}\right)+\frac{1}{a_np_n}.
\] 
But
\[\frac{1}{a_kp_k}-\frac{p_{k+1}}{a_{k+1}}=\frac{1}{a_kp_k}-\frac{p_{k+1}}{a_kp_kp_{k+1}+q_{k+1}}>\frac{1}{a_kp_k}-\frac{p_{k+1}}{a_kp_kp_{k+1}}=0.\]
\end{proof}
We end up the chapter with the simplest example of multibranched singularity, i.e. with a singularity defined locally by $x^d-y^d=0$ with $d\ge 2$. Its
link at singularity is the torus link $T_{d,d}$. Let us consider a set
\[\Sigma_d=\{\frac{i}{d}+\frac{j}{d},\,1\le i,j\le d-1\}.\]
Here the element $k/d$ appears in $\Sigma_d$ precisely $d-1-|d-1-k|$ times, according to possible presentations $k=i+j$, $1\le i,j\le d-1$.
Let $s_d(x)$ be the function computing the elements of $\Sigma_d$ in $(x,x+1)$ with a '$-$' sign and the others with '$+$' sign. Then
$s_d$ is almost equal to the \TLS{} of link $T_{d,d}$. We have the formula
\[s_d=2\sum_{k<d/2}(k-1)\chf{\frac{k}{d}}{\frac{d-k}{d}}-2(k-1)\sum_{k>d/2}\chf{\frac{d-k}{d}}{\frac{k}{d}}-(d-1).\]
The final term, $-(d-1)$, comes from the $d-1$ elements of the set $\Sigma_d$ of type $d/d$. They belong to any interval $(x,x+1)$.
Thus, the integral of $s_d$ is equal to
\[\int_0^1s_d=-2\sum_{k=1}^{d-1}(k-1)\frac{2k-d}{d}-(d-1).\]
But an elementary calculus shows that
\[
\sum_{k=1}^{d-1}(k-1)(2k-d)=
\frac{d(d-1)(d-2)}{6}.
\]
Hence
\[
\int_0^1s_d=-\frac13(d-1)(d+1).
\]
On the other hand, in order to resolve the singularity of $C$ we need only one blow-up. The exceptional divisor $E$ consists of single
rational curve with $E^2=-1$. Then $K=K_X=\alpha E$ and $C'=\beta E$ (as $E$ spans second (co)homology of blown-up space) and $K(K+E)=-2$
by genus formula, so $K=E$ and $C'\cdot E=d$, so $C'=-d\cdot E$. Thus $K+D=K+C'+E=(2-d)E$. Moreover, this divisor is nef, so its Zariski--Fujita
decomposition is trivial, $H=(2-d)E$, $N=0$,
so in this case
\[H^2=-(d-2)^2.\]
This shows that, in case of general links, a trivial analogue of Proposition~\ref{bound} does not hold. 
\begin{acknowledgement}
The author is very grateful to Tim Cochran and Stefan Friedl for explaining the rudiments of \rh invariants. He wishes also to
express his thanks to Andr\`as N\'emethi for various discussions on the subject.
\end{acknowledgement}

\end{document}